\documentstyle[amssymb,amscd,12pt]{amsart}

\newtheorem{lem}{Lemma}
\newtheorem{prop}{Proposition}

\theoremstyle{definition}

\theoremstyle{definition}
\newtheorem{thm}{Theorem}

\theoremstyle{remark}
\newtheorem{rem}{Remark}

\begin{document}

\newcommand{\thmref}[1]{Theorem~\ref{#1}}
\newcommand{\secref}[1]{Sect.~\ref{#1}}
\newcommand{\lemref}[1]{Lemma~\ref{#1}}
\newcommand{\propref}[1]{Proposition~\ref{#1}}
\newcommand{\corref}[1]{Corollary~\ref{#1}}
\newcommand{\remref}[1]{Remark~\ref{#1}}
\newcommand{\nc}{\newcommand}
\nc{\on}{\operatorname}
\nc{\ch}{\mbox{ch}}
\nc{\Z}{{\Bbb Z}}
\nc{\C}{{\Bbb C}}
\nc{\cond}{|\,}
\nc{\bib}{\bibitem}
\nc{\pone}{\Pro^1}
\nc{\pa}{\partial}
\nc{\F}{{\cal F}}
\nc{\arr}{\rightarrow}
\nc{\larr}{\longrightarrow}
\nc{\al}{\alpha}
\nc{\ri}{\rangle}
\nc{\lef}{\langle}
\nc{\W}{{\cal W}}
\nc{\gam}{\ovl{\gamma}}
\nc{\Q}{\ovl{Q}}
\nc{\q}{\widetilde{Q}}
\nc{\la}{\lambda}
\nc{\ep}{\epsilon}
\nc{\su}{\widehat{\frak{sl}}_2}
\nc{\gb}{\ovl{\frak{g}}}
\nc{\g}{\frak{g}}
\nc{\hh}{\ovl{\frak{h}}}
\nc{\h}{\frak{h}}
\nc{\n}{\frak{n}}
\nc{\ab}{\frak{a}}
\nc{\f}{\widehat{{\cal F}}}
\nc{\is}{{\bold i}}
\nc{\V}{{\cal V}}
\nc{\M}{\widetilde{M}}
\nc{\js}{{\bold j}}
\nc{\bi}{\bibitem}
\nc{\laa}{\ovl{\lambda}}
\nc{\fl}{B_-\backslash G}
\nc{\De}{\rho^\vee}
\nc{\G}{\widetilde{\frak{g}}}
\nc{\Li}{{\cal L}}
\nc{\fp}[2]{\frac{\pa}{\pa u_{#1}^{(#2)}}}
\nc{\Ve}{{\cal Vect}}
\nc{\sw}{\frak{sl}}
\nc{\La}{\Lambda}
\nc{\ds}{\displaystyle}
\nc{\HH}{\widetilde{\h}}
\nc{\bb}{\frak{b}}
\nc{\wt}{\widetilde}
\nc{\xx}{x}
\nc{\ovl}{\overline}

\title[Equivalence of two approaches to the mKdV
hierarchies]{Equivalence of two approaches to the mKdV hierarchies}

\author{Benjamin Enriquez}
\address{Centre de Math\'ematiques, URA 169 du CNRS, Ecole Polytechnique,
91128 \\ Palai\-seau, France}

\author{Edward Frenkel}
\address{Department of Mathematics, Harvard University,
Cambridge, MA 02138, USA}

\date{June 1996}

\maketitle

\begin{abstract}
The equivalence between the approaches of Drinfeld-Sokolov and
Fei\-gin-Frenkel to the mKdV hierarchies is established. A new derivation
of the mKdV equations in the zero curvature form is given. Connection with
the Baker-Akhiezer function and the tau-function is also discussed.
\end{abstract}

\section{Introduction.}
To each affine Kac-Moody algebra $\g$ one can associate a modified
Korteweg-de Vries (mKdV) hierarchy of non-linear partial differential
equations. The mKdV hierarchy, which can be viewed as a refined form of a
generalized KdV hierarchy (see \cite{DS}), is a completely integrable
hamiltonian system. The equations of the hierarchy can be written in
hamiltonian form, and the corresponding hamiltonian flows commute with each
other.

It is known that the equations of an mKdV hierarchy can be represented
in the zero curvature form
\begin{equation}    \label{zc}
[\pa_{t_n} + V_n,\pa_z + V] = 0,
\end{equation}
where $t_n$'s are the times of the hierarchy, and $t_1=z$. Here $V$ and
$V_n$ are certain time dependent elements of the centerless affine algebra
$\g$. To write $V$ explicitly, consider the principal abelian subalgebra
$\ab$ of $\g$ (the precise definition is given below). It has a basis $p_i,
i\in\pm I$, $I$ being the set of all exponents of $\g$ modulo the Coxeter
number. Then
$$V = p_{-1} + {\bold u}(z),$$ where ${\bold u}(z)$ lies in the Cartan
subalgebra $\h$ of $\g$.

The element $p_{-1}$ has degree $-1$ with respect to the principal
gradation of $\g$, while ${\bold u}$ has degree $0$. This makes finding an
element $V_n$ that satisfies \eqref{zc} a non-trivial problem. Indeed,
equation \eqref{zc} can be written as
\begin{equation}    \label{next}
\pa_{t_n} {\bold u} = [\pa_z + p_{-1} + {\bold u},V_n].
\end{equation}
The left hand side of \eqref{next} has degree $0$. Therefore $V_n$ should
be such that the expression in the right hand side of \eqref{next} is
concentrated in degree $0$.

Such elements can be constructed by the following trick (see
\cite{ZS,DS,W}). Suppose we found some ${\cal V}_n \in \g$ which satisfies
\begin{equation}    \label{next1}
[\pa_z + p_{-1} + {\bold u},{\cal V}] = 0.
\end{equation}
We can split ${\cal V}_n$ into the sum ${\cal V}_{+} + {\cal
V}_{-}$ of its components of positive and non-positive degrees with
respect to the principal gradation. Then $V_n = {\cal V}_{-}$ has the
property that the right hand side of \eqref{next} has degree $0$. Indeed,
from \eqref{next1} we find
$$[\pa_z + p_{-1} + {\bold u}(z),{\cal V}_{-}] = - [\pa_z + p_{-1} +
{\bold u}(z),{\cal V}_{+}],$$ which means that both commutators have
neither positive nor negative homogeneous components. Therefore
equation \eqref{next} makes sense. Now we have to find solutions of
equation \eqref{next1}.

Drinfeld and Sokolov \cite{DS} proposed a powerful method of finding
solutions of \eqref{next1}, which is closely related to the dressing method
of Zakharov and Shabat \cite{ZS}. Another approach was proposed by Wilson
\cite{W} (see also \cite{KW}).

Let us briefly explain the Drinfeld-Sokolov method. Recall that $\g$ has
the decomposition $\g = \n_+ \oplus \bb_-$, where $\n_+$ is the nilpotent
subalgebra of $\g$. Let $N_+$ be the corresponding Lie group. In \cite{DS}
it was proved that there exists an $N_+$--valued function $M(z)$, which is
called the dressing operator, such that
$$M(z)^{-1} \left( \pa_z + p_{-1} + {\bold u}(z) \right) M(z) = \pa_z +
p_{-1} + \sum_{i\in I} h_i(z) p_i,$$ where $h_i$'s are certain functions.

The dressing operator $M(z)$ is defined up to right multiplication by a
$z$--dependent element of the subgroup $A_+$ of $N_+$ corresponding to the
Lie algebra $\ab_+ = \ab \cap \n_+$. Thus, $M(z)$ represents a coset in
$N_+/A_+$. The element ${\cal V}_n = M(z) p_{-n} M(z)^{-1}$ clearly
satisfies \eqref{next1} and by substituting $V_n = (M(z) p_{-n}
M(z)^{-1})_-$ in equation \eqref{zc} for $n \in I$ one obtains the mKdV
hierarchy.

Recently, another approach to mKdV hierarchies was proposed by Feigin
and one of the authors \cite{FF:laws,FF}. In this approach, the flows
of the mKdV hierarchy are considered as vector fields on the space of
jets of the function ${\bold u}(z)$. Let
$\pi_0=\C[u_i^{(n)}]_{i=1,\ldots,l;n\geq 0}$, where $u_i=(\al_i,{\bold
u})$ and $u_i^{n} = \pa_z^n u_i$, be the ring of differential
polynomials in $u_i$'s. In \cite{FF}, $\pi_0$ was identified with the
ring of algebraic functions on the homogeneous space $N_+/A_+$. Thus,
each function ${\bold u}(z)$ gives rise to a function $K(z)$ with
values in $N_+/A_+$. The Lie algebra $\ab_- = \ab \cap \bb_-$
naturally acts on $N_+/A_+$ from the right. Consider the derivation
$\pa_n$ on $\pi_0$ which corresponds to the infinitesimal action of
$p_{-n}$ on $N_+/A_+$. These derivations clearly commute with each
other. Moreover, it was shown in \cite{FF} that $\pa_1$ coincides with
$\pa_z$ and therefore $\pa_n$'s are evolutionary (i.e. commuting with
$\pa_z$) derivations.

In this work we prove that the cosets $M(z)$ and $K(z)$, obtained by
the constructions of \cite{DS} and \cite{FF}, coincide. We then show
that the derivation $\pa_n$ satisfies equation \eqref{zc} with $V_n =
(K(z) p_{-n} K(z)^{-1})_-$. Thus, we establish an equivalence between
the two constructions. Note that another approach to establishing this
equivalence in the case of $\su$ based on KdV gauge fixing \cite{DS}
was proposed by one of the authors in \cite{E}.

This gives us a direct identification of the flows corresponding to $\pa_n$
and $\pa_{t_n}$.\footnote{In \cite{FF} the following indirect proof of this
fact was given: the derivations $\pa_n$ were identified in \cite{FF} with
the symmetries of the affine Toda equation corresponding to $\g$. But it is
known that mKdV equations constitute all symmetries of the affine Toda
equation, see \cite{DS,KW,W}.} Thus we obtain a new derivation of the zero
curvature representation of the mKdV hierarchies.

We remark that there exist generalizations of the mKdV hierarchies
which are associated to abelian subalgebras of $\g$ other then
$\ab$. It is known that the Drinfeld-Sokolov approach can be applied
to these generalized hierarchies \cite{DHM,HM}. On the other hand, the
approach of \cite{FF} can also be applied; in the case of the
non-linear Schr\"{o}dinger hierarchy, which corresponds to the
homogeneous abelian subalgebra of $\g$, this has been done by Feigin
and one of the authors \cite{FF:nls}. The results of our paper can be
extended to establish the equivalence between the two approaches in
this general context.

The paper is arranged as follows. In Sect.~2 we recall the
construction of \cite{FF} and derive the zero curvature equations. In
Sect.~3 we prove that the cosets $M$ and $K$ coincide and that the
derivations $\pa_n$ and $\pa_{t_n}$ coincide. We also explain the
connection with the KdV hierarchies. In Sect.~4 we construct a natural
system of coordinates on the group $N_+$ and using it give another
proof of the equivalence of two formalisms. Finally, in Sect.~5 we
obtain explicit formulas for the one-cocycles defined in \cite{FF} and
the densities of the hamiltonians of the mKdV hierarchy. We also
discuss a connection between the formalism of \cite{FF} and
$\tau$--functions.

\section{Unipotent cosets.}
\subsection{Notation.}
Let $\G$ be an affine algebra. It has generators $e_i,f_i,\al_i^\vee,
i=0,\ldots,l$, and $d$, which satisfy the standard relations
\cite{Kac}. The Lie algebra $\G$ carries a non-degenerate invariant inner
product $(\cdot,\cdot)$. One associates to $\G$ the labels $a_i, a_i^\vee,
i=0,\ldots,l$, the exponents $d_i, i=1,\ldots,l$, and the Coxeter number
$h$, see \cite{Kac}. We denote by $I$ the set of all positive integers,
which are congruent to the exponents of $\G$ modulo $h$ (with
multiplicities). The elements $e_i, i=0,\ldots,l$, and $f_i, i=0,\ldots,l$,
generate the nilpotent subalgebras $\n_+$ and $\n_-$ of $\G$,
respectively. The elements $\al_i^\vee$ generate the Cartan subalgebra
$\HH$ of $\G$. We have: $\G=\n_+ \oplus \bb_-$, where $\bb_- = \n_- \oplus
\HH$. Each $x \in \G$ can be uniquely written as $x_+ + x_-$, where $x_+
\in \n_+$ and $x_- \in \bb_-$.

The element $\ds C = \sum_{i=0}^l a_i^\vee \al_i^\vee$ of $\HH$ is a
central element of $\G$. Let $\g$ be the quotient of $[\G,\G]$ by $\C
C$. We identify $\G$ with the direct sum $\g \oplus \C C \oplus \C d$. The
Lie algebra $\g$ has a Cartan decomposition $\g = \n_+ \oplus \h \oplus
\n_-$, where $\h$ is spanned by $\al_i^\vee, i=1,\ldots,l$.

Set $\ds p_1 = \sum_{i=0}^l a_i e_i$. Let $\ab$ be the centralizer of $p_1$
in $\g$. This is an abelian subalgebra of $\g$ which we call the principal
abelian subalgebra. We have a decomposition: $\ab = \ab_+ \oplus \ab_-$,
where $\ab_+ = \ab \cap \n_+$, and $\ab_- = \ab \cap \bb_-$. It is known that
$\ab_\pm$ is spanned by elements $p_i, i\in \pm I$, which have degrees
$\deg p_i=i$ with respect to the principal gradation of $\g$. In
particular, we can choose $\ds p_{-1} = \sum_{i=0}^l
\frac{(\al_i,\al_i)}{2} f_i$, where $\al_i$'s are the simple roots of $\g$,
considered as elements of $\h$ using the inner product.

\begin{rem}    \label{multi}
For all affine algebras except $D_{2n}^{(1)}$, each exponent occurs exactly
once. In the case of $D_{2n}^{(1)}$, the exponent $2n-1$ has multiplicity
$2$. In this case, there are two generators of $\ab$, $p_i^1$ and $p_i^2$,
for $i$ congruent to $2n-1$ modulo the Coxeter number $4n-2$.\qed
\end{rem}

Let $N_+$ be the Lie group of $\n_+$. This is a prounipotent proalgebraic
group (see, e.g., \cite{FF}). The exponential map $\exp: \n_+ \arr N_+$ is
an isomorphism of proalgebraic varieties. Let $A_+$ be the image of $\ab_+$
under this map. The Lie algebra $\g$ acts on $N_+$ from the right because
$N_+$ can be embedded as an open subset in the flag manifold $\fl$ of
$\g$. Therefore the normalizer of $\ab_+$ in $\g$ acts on $N_+/A_+$ from
the right. In particular, $\ab_-$ acts on $N_+/A_+$, and each $p_{-n}, n\in
I$ gives rise to a derivation of $\C[N_+/A_+]$, see \cite{FF}. These
derivations commute with each other.

\subsection{Actions on the space of jets.}    \label{jets}
Consider the space $U$ of jets of a smooth function ${\bold u}(z): {\Bbb
A}^1 \arr \h$. The space $U$ is the inverse limit of the finite-dimensional
vector spaces $U_N = \on{span} \{ u_i^{(n)} \}_{i=1,\ldots,l;
n=1,\ldots,N}$, where $u_i = (\al_i,{\bold u})$, and $u_i^{(n)}=\pa_z^n
u_i$. Thus, the ring $\pi_0$ of regular functions on $U$ is
$\C[u_i^{(n)}]_{i=1,\ldots,l; n\geq 0}$. The derivative $\pa_z$ gives rise
to a derivation of $\pi_0$.

\begin{thm}[\cite{FF}, Theorem 2]    \label{iso}
{\em There is an isomorphism of rings $$\C[N_+/A_+] \simeq \C[u_i^{(n)}],$$
under which $p_{-1}$ gets identified with $\pa_z$.}
\end{thm}

Let $\pa_n$ be the derivation of $\C[u_i^{(n)}]$ corresponding to
$p_{-n}$ under this isomorphism. The theorem shows that the
derivations $\pa_n$ are evolutionary, i.e. commuting with $\pa_z$. We
would like to represent the action of these derivations on ${\bold
u}(z)$ explicitly in the zero curvature form \eqref{zc}.

For $g \in G$ and $x \in \g$ we will write $g x g^{-1}$ for $\on{Ad}_g(x)$.

\begin{prop}    \label{zcgen}
For $K \in N_+/A_+$,
\begin{equation}    \label{zcgenf}
[\pa_m + (Kp_{-m}K^{-1})_-,\pa_n + (Kp_{-n}K^{-1})_-] = 0, \quad \quad
\forall m,n \in I.
\end{equation}
\end{prop}

Let us explain the meaning of formula \eqref{zcgenf}. For each $K \in
N_+/A_+$, $K p_{-n}K^{-1}$ is a well-defined element of $\g$. The Lie
algebra $\g$ can be realized as $\ovl{\g}\otimes \C((t))$ (or a subalgebra
thereof if $\g$ is twisted) for an appropriate finite-dimensional Lie
algebra $\ovl{\g}$. If we choose a basis in $\ovl{\g}$, we can consider an
element of $\g$ as a set of Laurent power series. In particular, for $K
p_{-n}K^{-1}$, any Fourier coefficient of each of these power series gives
us an algebraic function on $N_+/A_+$. Hence, by \thmref{iso}, each
coefficient corresponds to a differential polynomial in $u_i$'s, and we can
apply $\pa_m$ to it.

In order to prove formula \eqref{zcgenf}, we need to find an explicit
formula for the action of $\pa_n$ on $Kp_{-m}K^{-1}$.

Let us first obtain a formula for the infinitesimal action of an element of
$\g$ on $N_+$. Recall from \cite{FF} that since $N_+$ embeds as an open
subset in the flag manifold $\fl$, the Lie algebra $\g$ infinitesimally
acts on $N_+$ from the right by vector fields. Therefore $\g$ acts by
derivations on $\C[N_+]$. Since $N_+$ acts on $\g$, we obtain a
homomorphism from $N_+$ to the group of automorphisms of $\ovl{\g}$ over
the ring $\C[[t]]$.

Now we can consider each element of $N_+$ as a matrix, whose entries are
Taylor power series. Each Fourier coefficient of such a series defines an
algebraic function on $N_+$, and the ring $\C[N_+]$ is generated by these
functions. Hence any derivation of $\C[N_+]$ is uniquely determined by its
action on these functions. We can write this action concisely as follows:
$\nu \cdot x = y$, where $x$ is a ``test'' matrix representing an element
of $N_+$, and $y$ is another matrix, whose entries are the results of the
action of $\nu$ on the entries of $x$.

For $a \in \g$, let $a^R$ be the derivation of $\C[N_+]$ corresponding to
the right infinitesimal action of $a$ on $N_+$. For $b \in \n_+$, let $b^L$
be the derivation of $\C[N_+]$ corresponding to the left infinitesimal
action of $b$ on $N_+$.

\begin{lem}
\begin{align}    \label{actionr}
a^R \cdot \xx &= (\xx a\xx^{-1})_+ \xx, \quad \quad \forall a \in \g,\\
\label{actionl} b^L \cdot \xx &= b\xx, \quad \quad \forall b \in \n_+.
\end{align}
\end{lem}

\begin{pf} Consider a one-parameter subgroup $a(\ep)$ of $G$, such that
$a(\ep) = 1 + \ep a + o(\ep)$. We have: $\xx \cdot a(\ep) = \xx + \ep \xx
a + o(\ep)$. For small $\ep$ we can factor $\xx \cdot a(\ep)$ into a
product $y_- y_+$, where $y_+ = \xx + \ep y_+^{(1)} + o(\ep) \in N_+$ and
$y_- = 1 + \ep y_-^{(1)} \in B_-$. We then find that $y_-^{(1)} \xx +
y_+^{(1)} = \xx a$, from which we conclude that $y_+^{(1)} = (\xx
a\xx^{-1})_+ \xx$. This proves formula \eqref{actionr}. Formula
\eqref{actionl} is obvious.
\end{pf}

It follows from formula \eqref{actionr} that
\begin{equation}    \label{actionc}
a^R \cdot \xx  v \xx^{-1} = [(\xx a\xx^{-1})_+,\xx  v \xx^{-1}], \quad
\quad a,v \in \g.
\end{equation}
If $a$ and $v$ are both elements of $\ab$, then formula \eqref{actionc}
does not change if we multiply $\xx$ from the right by an element of
$A_+$. Denote by $K$ the coset of $\xx$ in $N_+/A_+$. Then we can write:
\begin{equation}    \label{actionc1}
\pa_n \cdot K v K^{-1} = [(Kp_{-n}K^{-1})_+,K v K^{-1}], \quad \quad v
\in \ab.
\end{equation}

\vspace*{5mm}
\noindent {\em Proof of \propref{zcgen}.} Substituting $v=p_{-m}$ into
formula \eqref{actionc1}, we obtain:
$$\pa_n \cdot K p_{-m} K^{-1} = [(Kp_{-n}K^{-1})_+,K p_{-m} K^{-1}].$$
Hence $$\pa_n \cdot (K p_{-m} K^{-1})_- = [(Kp_{-n}K^{-1})_+,K p_{-m}
K^{-1}]_- = [(Kp_{-n}K^{-1})_+,(K p_{-m} K^{-1})_-]_-.$$ Therefore we
obtain:
\begin{align*}
& [\pa_m + (Kp_{-m}K^{-1})_-,\pa_n + (Kp_{-n}K^{-1})_-] \\ = & \pa_m \cdot
(K p_{-n}K^{-1})_- - \pa_n \cdot (K p_{-n} K^{-1})_- + [(K p_{-m}
K^{-1})_-,(K p_{-n} K^{-1})_-] \\ = & [(Kp_{-m}K^{-1})_+,(K p_{-n}
K^{-1})_-]_- - [(Kp_{-n}K^{-1})_+,K p_{-m} K^{-1}]_- \\ + & [(K p_{-m}
K^{-1})_-,(K p_{-n} K^{-1})_-]_-.
\end{align*}

Adding up the first and the last terms, we obtain $$[K p_{-m}
K^{-1},(Kp_{-n}K^{-1})_-]_- - [(Kp_{-n}K^{-1})_+,K p_{-m} K^{-1}]_- =$$
$$= [K p_{-m} K^{-1},K p_{-n} K^{-1}]_- = 0,$$
and the proposition is proved.\qed

\subsection{Zero curvature form.}
In order to write the action of $\pa_n$ in the zero curvature form, we will
use \propref{zcgen} in the case $m=1$. But first we determine
$(Kp_{-1}K^{-1})_-$ explicitly.

\begin{lem}    \label{p-1}
\begin{equation}    \label{boldu}
(Kp_{-1}K^{-1})_- = p_{-1} + {\bold u}.
\end{equation}
\end{lem}

It is clear that $(Kp_{-1}K^{-1})_- = p_{-1} + x$, where $x \in \h$. Hence
we need to show that $x={\bold u}$, or, equivalently, that $(\al_i,x) =
u_i, i=1,\ldots,l$. We can rewrite the latter formula as
$u_i=(\al_i,Kp_{-1}K^{-1})_-$, and hence as
$u_i=(\al_i,Kp_{-1}K^{-1})$. To establish the last formula, recall the
interpretation of $u_i$ from \cite{FF}.

Consider the module $M_\la^*$ contragradient to the Verma module
$M_\la$ over $\g$ with highest weight $\la$. This module can be
realized in the space $\C[N_+]$ in such a way that the highest weight
vector $v_\la$ corresponds to the constant function. For $a \in \g$
denote by $f_\la(a)$ the function on $N_+$ which corresponds to $a
\cdot v_\la$. Then $u_i = f_{\al_i}(p_{-1})$ \cite{FF}. But in fact
there is a general formula for $f_\la(a)$ due to Kostant \cite{Kos}.

\begin{prop}[\cite{Kos}, Theorem 2.2]    \label{kost}
Consider $\la \in \h^*$ as a functional on $\g$ which is trivial on
$\n_\pm$. Let $\langle \cdot,\cdot \rangle$ be the pairing between $\g^*$
and $\g$. Then $f_\la(a)(\xx) = \langle \la,\xx a \xx^{-1} \rangle$.
\end{prop}

The formula above immediately implies that the function $u_i$ on $N_+/A_+$
takes value $(\al_i,Kp_{-1}K^{-1})$ at $K \in N_+/A_+$. This completes the
proof of \lemref{p-1}.

Now specializing $m=1$ in formula \eqref{zcgenf} and using \lemref{p-1} we
obtain the zero curvature representation of the equations.

\begin{thm}
\begin{equation}    \label{dsmkdv1}
[\pa_z + p_{-1} + {\bold u},\pa_n + (Kp_{-n}K^{-1})_-] = 0.
\end{equation}
\end{thm}

This equation can be rewritten as
\begin{equation}    \label{dsmkdv}
\pa_n {\bold u} = \pa_z (Kp_{-n}K^{-1})_- + [p_{-1}+{\bold
u},(Kp_{-n}K^{-1})_-].
\end{equation}

The map $K \arr K p_{-n}K^{-1}$ defines an embedding of $N_+/A_+$ into $\g$
as an $N_+$--orbit. The entries of the matrix $Kp_{-n}K^{-1}$ are Laurent
series in $t$ whose coefficients are differential polynomials in $u_i,
i=1,\ldots,l$ (see the paragraph after \propref{zcgen}).  Equation
\eqref{dsmkdv} expresses $\pa_n u_i$ in terms of differential polynomials
in $u_i$'s. Since, by construction, $\pa_n$ commutes with $\pa_1 \equiv
\pa_z$, formula \eqref{dsmkdv} uniquely determines $\pa_n$ as an
evolutionary derivation of $\C[u_i^{(n)}]$.

\section{Equivalence with the Drinfeld-Sokolov formalism.}

\subsection{Identification of cosets.}
For $v \in \ab$ and $K \in N_+/A_+$ set
\begin{equation}    \label{calV}
{\cal V}_v = K v K^{-1}.
\end{equation}
Since $\ab$ is commutative, this is a well-defined element of $\g$.

\begin{prop}    \label{centr}
\begin{equation}    \label{comm}
[\pa_z + p_{-1} + {\bold u},{\cal V}_v]=0, \quad \quad \forall v \in
\ab.
\end{equation}
\end{prop}

\begin{pf} Using formula \eqref{actionc1} and \lemref{p-1} we obtain:
\begin{equation}    \label{dv}
\pa_z {\cal V}_v = [(Kp_{-1}K^{-1})_+,{\cal V}_v] =  -
[(Kp_{-1}K^{-1})_-,{\cal V}_v] = - [p_{-1} + {\bold u},{\cal V}_v].
\end{equation}
\end{pf}

Now we define the Drinfeld-Sokolov dressing operator $M$.

\begin{prop}[\cite{DS}, Proposition 6.2]    \label{conj}
There exists an element $M=M(z) \in N_+$, such that
\begin{equation}    \label{ds}
M^{-1} \left( \pa_z + p_{-1} + {\bold u}(z) \right) M = \pa_z + p_{-1}
+ \sum_{i\in I} h_i p_i,
\end{equation}
where $h_i$'s are functions. $M$ is defined uniquely up to right
multiplication by a (possibly $z$--dependent) element of $A_+$. One can
choose $M$ in such a way that all entries of its matrix and all $h_i$'s are
polynomials in $u_i^{(n)}, i=1,\ldots,l; n\geq 0$.
\end{prop}

The proposition defines a map from the space of smooth functions
${\bold u}(z): {\Bbb A}^1 \arr \h$ to the space of smooth functions
${\Bbb A}^1 \arr N_+/A_+$, ${\bold u}(z) \arr M(z)$. On the other
hand, \thmref{iso} also defines such a map ${\bold u}(z) \arr K(z)$. The
following lemma will allow us to identify these two maps.

\begin{rem} Note that both maps are local in the following sense. For
each $z$, $M(z)$ and $K(z)$ depend only on the jet of ${\bold u}$ at
$z$. In particular, for each $v \in \ab$, all entries of the matrices
$M(z)vM(z)^{-1}$ and $K(z)vK(z)^{-1}$ are Taylor series whose
coefficients are differential polynomials in $u_i$'s.\qed
\end{rem}

\begin{lem}[\cite{DS}]    \label{unique}
Let ${\cal V}$ be an element of $\g$ of the form ${\cal V} = p_{-n} +$
terms of degree higher than $-n$ with respect to the principal
gradation on $\g$, such that
\begin{equation}    \label{nov}
[\pa_z + p_{-1} + {\bold u},{\cal V}] = 0.
\end{equation}
Then ${\cal V} = M v M^{-1}$, where $M \in N_+$ satisfies
\eqref{ds} and $v \in \ab$ is such that $v = p_{-n} +$ terms of
degree higher than $-n$.
\end{lem}

The proof of the lemma requires the following important result.

\begin{prop}[\cite{Kac1}, Proposition 3.8]    \label{kac}
The Lie algebra $\g$ has the decomposition $\g = \ab \oplus \on{Im}
(\on{ad} p_{-n})$ for each $n \in I$. Moreover, $\on{Ker} (\on{ad} p_{-n})
= \ab$.
\end{prop}

\noindent {\em Proof of \lemref{unique}}. If ${\cal V}$ satisfies
\eqref{nov}, then we obtain from \propref{conj}:
\begin{equation}    \label{nov1}
[\pa_z + p_{-1} + \sum_{i\in I} h_i p_i,M^{-1}{\cal V}M] = 0.
\end{equation}
We can write $M^{-1}{\cal V}M$ as a sum $\sum_j v_j$ of its homogeneous
components of principal degree $j$. According to \propref{kac}, each $v_j$
can be split into the sum of $v_j^0 \in \ab_+$ and $v_j^1 \in \on{Im}
(\on{ad} p_{-1})$.

Suppose that $M^{-1}{\cal V}M$ does not lie in $\ab_+$. Let $j_0$ be the
smallest number such that $v_{j_0}^1 \neq 0$. Then the term of smallest
degree in \eqref{nov1} is $[p_{-1},v_{j_0}^1]$ which is non-zero, because
$\on{Ker} (\on{ad} p_{-n}) = \ab_+$. Hence \eqref{nov1} can not hold.

Therefore $M^{-1}{\cal V}M \in \ab$. But then \eqref{nov1} gives: $\pa_z
v_j = 0$ for all $j$. This means that each $v_j$ is a constant element of
$\ab$, and Lemma follows.\qed

\begin{thm}    \label{ident}
{\em The cosets $M(z)$ and $K(z)$ in $N_+/A_+$ assigned in \cite{DS}
and \cite{FF}, respectively, to the jet of function ${\bold u}: {\Bbb
A}^1 \arr \h$ at $z$, coincide.}
\end{thm}

\begin{pf}
According to \propref{centr}, $$[\pa_z + p_{-1} + {\bold u},Kp_{-n}K^{-1}]
=0.$$ Since $Kp_{-n}K^{-1} = p_{-n} +$ terms of degree higher than $-n$
with respect to the principal gradation, we obtain from \lemref{unique}
that $Kp_{-n}K^{-1} = M v M^{-1}$, where $M \in N_+$ satisfies \eqref{ds}
and $v \in \ab$. This implies that $v=p_{-n}$ and that $M=K$ in $N_+/A_+$.

Indeed, from the equality $$Kp_{-n}K^{-1} = M v M^{-1}$$ we obtain that
$(M^{-1} K) p_{-n} (M^{-1} K)^{-1}$ lies in $\ab$. We can represent $M^{-1}
K$ as $\exp y$ for some $y \in \n_+$. Then $(M^{-1} K) p_{-n} (M^{-1}
K)^{-1} = v$ can be expressed as a linear combination of multiple
commutators of $y$ and $p_{-n}$: $$e^y p_{-n} (e^y)^{-1} = \sum_{n\geq 0}
\frac{1}{n!}  (\on{ad} y)^n \cdot p_{-n}.$$ We can write $y = \sum_{j>0}
y_j$, where $y_j$ is the homogeneous component of $y$ of principal degree
$j$. It follows from \propref{kac} that $\n_+ = \ab_+ \oplus \on{Im}
(\on{ad} p_{-n})$. Therefore each $y_j$ can be further split into a sum of
$y_j^0 \in \ab_+$ and $y_j^1 \in \on{Im} (\on{ad} p_{-1})$.

Suppose that $y$ does not lie in $\ab_+$. Let $j_0$ be the smallest number
such that $y_{j_0}^1 \neq 0$. Then the term of smallest degree in $e^y
p_{-n} (e^y)^{-1}$ is $[y_{j_0}^1,p_{-n}]$ which lies in $\on{Im} (\on{ad}
p_{-n})$ and is non-zero, because $\on{Ker} (\on{ad} p_{-n}) =
\ab_+$. Hence $e^y p_{-n} (e^y)^{-1}$ can not be an element of $\ab_+$.

Therefore $y \in \ab_+$ and so $M^{-1} K \in A_+$, which means that $M$ and
$K$ represent the same coset in $N_+/A_+$, and that $v=p_{-n}$.
\end{pf}

\subsection{Identification of the equations.}
As was explained in the previous section, \thmref{iso} allows us to define
a set of commuting derivations $\pa_n, n\in I$, of $\pi_0$, or
equivalently, vector fields on the space of jets $U$. These derivations can
be represented in the zero curvature form \eqref{dsmkdv1}.

On the other hand, in \cite{DS} another set of derivations $\pa_{t_n}, n\in
I$, of $\pi_0$ was defined in the zero curvature form. Set
\begin{equation}    \label{vn}
V_n = \left( M(z) p_{-n} M(z)^{-1} \right)_-,
\end{equation}
where $M(z)$ is defined in \propref{conj}. In particular,
formula \eqref{ds} shows that $V_1 = V = p_{-1} + {\bold u}$. The $n$th
zero curvature equation is equation \eqref{zc}. Now we obtain from
\thmref{ident}

\begin{thm}    \label{dsff}
{\em The derivations $\pa_n$ and $\pa_{t_n}$ coincide.}
\end{thm}

\begin{rem}    \label{intcurves}
This theorem together with \thmref{iso} implies that solutions of the
mKdV hierarchy are just the integral curves of the vector fields of
the infinitesimal action of the Lie algebra $\ab_-$ on $N_+/A_+$.\qed
\end{rem}

\begin{rem} The variable $t^{-1}$ appearing in the affine algebra $\g$ is
often denoted by $\la$, and is called the spectral parameter.\qed
\end{rem}

\begin{rem}
For each $n \in I \cup -I$, the map $K \mapsto K p_{-n} K^{-1}$ defines an
embedding $N_+/A_+ \arr \g$, because the stabilizer of $p_{-n}$ in $N_+$ is
$A_+$. In practice, it is convenient to find $K p_{-n} K^{-1}$ using
equation
\begin{equation}    \label{recu}
[\pa_z + p_{-1} + {\bold u}(z),Kp_{-n}K^{-1}]=0,
\end{equation}
which follows from formula \eqref{comm}. We can split $Kp_{-n}K^{-1}$ into
the sum of homogeneous components lying in $\ab$ and in $\on{Im} (\on{ad}
p_{-1})$. These homogeneous components can then be determined recursively
using equation \eqref{recu} as explained in \cite{W}, Sect.~3.

This recursion is actually non-trivial: at certain steps one has to take
the anti-derivative of a differential polynomial. But we know from
\propref{centr} that the element $Kp_{-n}K^{-1}$ satisfies \eqref{recu} and
that its entries are differential polynomials (see the paragraph after
\propref{zcgen}). Therefore whenever an anti-derivative occurs, it can be
resolved in the ring of differential polynomials. Another proof of this
fact has been given by Wilson \cite{W}.

Every time we compute the anti-derivative, we have the freedom of adding an
arbitrary constant. This corresponds to adding to $Kp_{-n}K^{-1}$ a linear
combination of $Kp_mK^{-1}$ with $m>-n$.\qed
\end{rem}

\begin{rem} The map which attaches to $u_i$'s a coset in $N_+/A_+$ can
be viewed as a universal feature in various approaches to soliton
equations. In this section we have explained how these maps arise in
the formalisms of \cite{DS} and \cite{FF} and proved that these maps
coincide.

But a map to $N_+/A_+$ can also be found, in a somewhat disguised form,
in the approach to the soliton equations based on Sato's Grassmannian,
see \cite{SW,Wil}. One can associate to $u_i$'s their Baker-Akhiezer
function $\Psi$ which is a solution of the equation
\begin{equation}    \label{baker}
(\pa_z + p_{-1} + {\bold u}(z)) \Psi = 0,
\end{equation}
and more generally the equations
\begin{equation}    \label{bakergen}
(\pa_n + (Kp_{-n}K^{-1})_-) \Psi = 0, \quad \quad \forall n \in I.
\end{equation}

In our notation, Segal and Wilson \cite{SW,Wil} attach in the case of
$\g=\sw_n$ a Baker-Akhiezer function $\Psi$ to each point $x$ of the
flag manifold $B_- \backslash G$ using its realization via an infinite
Grassmannian. The flows of the mKdV hierarchy then correspond to the
right infinitesimal action of $\ab_-$ on the flag manifold. As $x$
moves along the integral curves of the vector fields $\pa_n$, the
Baker-Akhiezer function evolves according to the mKdV hierarchy and so
does the function ${\bold u}$. One shows \cite{SW} that $\Psi$ is
regular at a given set of times of the hierarchy if the corresponding
point of the flag manifold lies in the big cell (which is isomorphic
to $N_+$). Moreover, ${\bold u}$ does not change under the right
action of $A_+$ on $x$ \cite{Wil}. Thus, one obtains a map which
assigns to ${\bold u}$ an element of $N_+/A_+$.

In \cite{Wil} the equivalence between the dressing method and the
Grassmannian approach was established (see also \cite{HM}). Therefore
this map coincides with the map studied in our paper.

We will derive an explicit formula for the Baker-Akhiezer function in
Sect.~4.\qed
\end{rem}

\subsection{From mKdV to KdV} 
First let us recall the definition of the KdV hierarchies from
\cite{DS}. Consider the operator
\begin{equation}    \label{kdvL}
\pa_z + p_{-1} + {\bold u}(z),
\end{equation}
where now ${\bold u}(z)$ lies in the finite-dimensional Borel
subalgebra $\h \oplus \ovl{\n}_+$, where $\ovl{\n}_+$ is generated by
$e_i, i=1,\ldots,l$. Drinfeld and Sokolov construct in \cite{DS} the
dressing operator and the zero-curvature equations \eqref{zc} for this
operator in the same way as for the mKdV hierarchy using formulas
\eqref{ds} and \eqref{vn}.

The Lie group $\ovl{N}_+$ of $\ovl{\n}_+$ acts naturally on the space
of operators \eqref{kdvL} and these equations preverve the
corresponding gauge equivalence classes \cite{DS}. Thus one obtains a
system of compatible evolutionary equations on the gauge equivalence
classes, which is called the generalized KdV hierarchy corresponding
to $\G$. Let $\ovl{\n}_+^0$ be a subspace of $\ovl{\n}_+$ that is
transversal to the image in $\ovl{\n}_+$ of the operator $\on{ad}
\ovl{p}_{-1}$, where $\ovl{p}_{-1} = \sum_{i=1}^l
\frac{(\al_i,\al_i)}{2} f_i$. It is shown in \cite{DS} that each
equivalence class contains a unique operator \eqref{kdvL} satisfying
the condition that ${\bold u} \in \ovl{\n}_+^0$.

The space $\ovl{\n}_+^0$ is $l$--dimensional. If we choose coordinates
$v_1,\ldots,v_l$ of $\ovl{\n}_+^0$, then the KdV equations can be
written as partial differential equations on $v_i$'s. On the other
hand, the dressing operator $M(z)$ corresponding to a gauge class of
operators \eqref{kdvL} should now be considered as a double coset in
$\ovl{N}_+ \backslash N_+/A_+$. Thus, a smooth function ${\bold v}(z)
= (v_1(z),\ldots,v_l(z)): {\Bbb A}^1 \arr \ovl{\n}_+^0$ gives rise to
a smooth function ${\Bbb A}^1 \arr \ovl{N}_+ \backslash N_+/A_+$.

Denote by ${\cal L}$ the space of all operators \eqref{kdvL} where
${\bold u} \in \h$, and by $\wt{{\cal L}}$ the space of all operators
\eqref{kdvL} where ${\bold u} \in \ovl{\n}_+^0$. We obtain a
surjective map ${\cal L} \arr \wt{{\cal L}}$, which sends an operator
from ${\cal L}$ to the unique representative of its gauge class lying
in $\wt{{\cal L}}$. This map is called the Miura transformation. It
induces a homomorphism of differential rings $\C[v_i^{(n)}] \arr
\C[u_i^{(n)}]$.

It was shown in \cite{FF:laws} that the image of $\C[v_i^{(n)}]$ in
$\C[u_i^{(n)}]$ coincides with the invariant subspace of
$\C[u_i^{(n)}]$ under the left action of the group $\ovl{N}_+$. Hence
we obtain from \thmref{iso} that $\C[v_i^{(n)}] \simeq \C[\ovl{N}_+
\backslash N_+/A_+]$. Thus, we obtain a local map which assigns to
each smooth function ${\bold v}(z): {\Bbb A}^1 \arr \ovl{\n}_+^0$ a
smooth function ${\Bbb A}^1 \arr \ovl{N}_+ \backslash
N_+/A_+$. According to the results of this section, this map coincides
with the Drinfeld-Sokolov map defined above. We also see that the KdV
flows on $\ovl{N}_+ \backslash N_+/A_+]$ correspond to the right
infinitesimal action of $\ab_-$ on it considered as an open subset of
the ``loop space'' $\ovl{G}[t^{-1}]\backslash G/A_+$. Thus, the
passage from mKdV hierarchy to the KdV hierarchy simply consists of
projecting from the flag manifold $B_- \backslash G$ to the loop space
$\ovl{G}[t^{-1}] \backslash G$.

\begin{rem}
Drinfeld and Sokolov attached in \cite{DS} a generalized KdV hierarchy
to each vertex of the Dynkin diagram of $\G$; the hierarchy considered
above corresponds to the $0$th node. In general, we obtain the
following picture.

Fix $j$ between $0$ and $l$. Let $\ovl{\n}_+^{j}$ be the
finite-dimensional Lie subalgebra of $\n_+$ generated by $e_i, i\neq
j$. Let $\ovl{N}_+^{j}$ be the corresponding Lie subgroup of
$N_+$. The dressing operator of the $j$th generalized KdV hierarchy
gives rise to a double coset in $\ovl{N}_+^{j} \backslash
N_+/A_+$. On the other hand, there is an isomorphism between
$\C[\ovl{N}_+^{j} \backslash N_+/A_+]$ and the ring of
$\ovl{\n}_+^{j}$--invariants of $\C[u_i^{(n)}]$ (with respect to the
left action). The latter is itself a ring of differential polynomials
in $l$ variables. Note that it coincides with the intersection of
kernels of the operators $e_i^L, i\neq j$, which are classical limits
of the so-called screening operators (see \cite{FF:laws}).\qed
\end{rem}

\section{Realization of $\C[N_+]$ as a polynomial ring.}
The approach to the mKdV and affine Toda equations used in \cite{FF} and
here is based on \thmref{iso} which identifies $\C[N_+/A_+]$ with the ring
of differential polynomials $\C[u_i^{(n)}]_{i=1,\ldots,l;n\geq 0}$. In this
section we add to the latter ring new variables corresponding to $A_+$
and show that the larger ring thus obtained is isomorphic to
$\C[N_+]$. An analogous construction has been given in \cite{EF} in the
lattice case.

\subsection{Coordinates on $N_+$.}
Consider $u_i^{(n)}, i=1,\ldots,l; n\geq 0$, as $A_+$--invariant
regular functions on $N_+$. Recall that $$u_i(\xx )=(\al_i,\xx
p_{-1}\xx^{-1}), \quad \quad \xx \in N_+,$$ Now choose an element
$\chi$ of $\HH$, such that $(\chi,C) \neq 0$. Introduce the regular
functions $\chi_n, n\in I$, on $N_+$ by the formula:
\begin{equation}    \label{chin}
\chi_n(\xx )=(\chi,\xx p_{-n}\xx^{-1}), \quad \quad \xx \in N_+.
\end{equation}
Note that here we consider the action of $N_+$ on $\G$ and the pairing
on $\G$.

\begin{thm}    \label{identif}
$\C[N_+] \simeq \C[u_i^{(n)}]_{i=1,\ldots,l;n\geq 0} \otimes
\C[\chi_n]_{n\in I}.$
\end{thm}

\begin{pf} Let us show that the functions $u_i^{(n)}$'s and $\chi_n$'s
are algebraically independent. In order to do that, let us compute the
values of the differentials of these functions at the origin. Those are
elements of the cotangent space to the origin, which is isomorphic to the
dual space $\n_+^*$ of $\n_+$.

It follows from \propref{kac} that $\n_+^*$ can be written as $\n_+^*
= \ab_+^* \oplus \wt{\n}_+^*$, where $\ab_+^* = \on{Ker} (\on{ad}^*
p_{-1})$ and $\wt{\n}_+^*$ is the annihilator of $\ab_+$ with respect
to the pairing between $\n_+$ and $\n_+^*$. Moreover, if we decompose
$\wt{\n}_+^*$ with respect to the principal gradation as
$\oplus_{j=1}^\infty \wt{\n}_+^{*,j}$, then $\dim \wt{\n}_+^{*,j}=l$
for all $j>0$, and $\on{ad}^* p_{-1}$ maps $\wt{\n}_+^{*,j}$
isomorphically to $\wt{\n}_+^{*,j-1}$ for $j>1$.

By construction of $u_i$'s given in \cite{FF}, $du_i|_1, i=1,\ldots,l$,
form a basis of $\wt{\n}_+^{*,1}$, and hence $du^{(n)}_i|_1,
i=1,\ldots,l$, form a basis of $\wt{\n}_+^{*,n}$. Thus, the covectors
$du^{(n)}_i|_1, i=1,\ldots,l; n\geq 0$, are linearly independent. Let us
show now that the covectors $d\chi_n|_1$ are linearly independent from them
and among themselves. For that it is sufficient to show that the pairing
between $dF_m|_1$ and $p_n$ is non-zero if and only if $n=m$. But we have:
\begin{equation}    \label{cox}
p_n^R \cdot (\chi,\xx p_{-m}\xx^{-1}) = (\chi,\xx[p_n,p_{-m}]\xx^{-1})
\end{equation}
$$= (\chi,n(p_n,p_{-n})C) \delta_{n,-m} = n(p_n,p_{-n})(\chi,C)
\delta_{n,-m},$$ where $h$ is the Coxeter number. Therefore this pairing
equals $n (p_n,p_{-n}) (\chi,C) \delta_{n,-m}$. This satisfies the
condition above.

Thus, the functions $u_i^{(n)}$'s and $\chi_n$'s are algebraically
independent. Hence we have an embedding $\C[u_i^{(n)}]_{i=1,\ldots,l;n\geq
0} \otimes \C[\chi_n]_{n\in I} \arr \C[N_+]$. But the characters of the two
spaces with respect to the principal gradation are both equal to
$$\prod_{n\geq 0} (1-q^n)^{-l} \prod_{i\in I} (1-q^i)^{-1}.$$ Hence this
embedding is an isomorphism.
\end{pf}

\subsection{Another proof of \thmref{ident}.}
Now we will explain another point of view on the equivalence between the
formalisms of \cite{DS} and \cite{FF} established in Sect.~3.

In any finite-dimensional representation of $N_+$, each element $x$ of
$N_+$ is represented by a matrix whose entries are Taylor series in $t$
with coefficients in $\C[N_+]$.

Denote $c_n = (n (p_n,p_{-n}) (\chi,C))^{-1}$.

\begin{prop}    \label{localm}
Let $x$ be an element of $N_+$. We associate to it another element of
$N_+$, $$\ovl{x} = x \exp \left( - \sum_{n\in I} c_n p_{n} \chi_n(x)
\right)$$ In any finite-dimensional representation of $N_+$, $\ovl x$ is
represented by a matrix whose entries are Taylor series with coefficients
in the ring of differential polynomials in $u_i, i=1,\ldots,l$

The map $N_+ \arr N_+$ which sends $x$ to $\ovl{x}$ is constant on the
right $A_+$--cosets, and hence defines a section $N_+/A_+ \arr N_+$.
\end{prop}

\begin{pf}
Each entry of $$\ovl{x} = x \exp \left( - \sum_{n\in I}
c_n p_{n} \chi_n(x) \right)$$ is a function on
$N_+$. According to \thmref{iso} and \thmref{identif}, to prove the
proposition it is sufficient to show that each entry of $\ovl{x}$ is
invariant under the right action of $\ab_+$. By formula \eqref{cox} we
obtain for each $m \in I$:
$$
p_m^R \cdot \chi_n = c_n^{-1} \delta_{n,-m}, 
$$
and hence
$$
p_m^R \left( x \exp \left( - \sum_{n\in I} c_n p_{n}\chi_n(x) \right)
\right) =$$ $$x p_m \exp \left( - \sum_{n\in I} c_n p_{n}\chi_n \right) + x
\exp \left( - \sum_{n\in I} c_n p_{n}\chi_n \right) (-p_m) = 0.
$$
Therefore $\ovl{x}$ is right $\ab_+$--invariant.

To prove the second statement, let $a$ be an element of $A_+$ and let
us show that $\overline{xa} = \ovl{x}$. We can write: $a = \exp
\left(\sum_{i \in I} \al_i p_i \right)$. Then according to formulas
\eqref{chin} and \eqref{cox}, $\chi_n(xa) = (\chi,xa p_{-n} a^{-1}
x^{-1}) = (\chi,x p_{-n} x^{-1}) + c_n^{-1} \al_n = \chi_n(x) +
c_n^{-1} \al_n$. Therefore
$$\overline{xa} = xa \exp \left( - \sum_{n \in I} \al_n p_n - \sum_{n
\in I} c_n p_n \chi_n(x) \right) = \ovl{x}.$$
\end{pf}

Consider now the matrix $\ovl{x}$. According to \propref{localm}, the
entries of $\ovl{x}$ are Taylor series with coefficients in differential
polynomials in $u_i$'s. Hence we can apply to $\ovl{x}$ any derivation of
$\C[u_i^{(n)}]$, in particular, $\pa_n = p_{-n}^R$. In the following
proposition we consider $p_i$ and ${\bold u}$ as matrices acting in a
finite-dimensional representation.

\begin{lem}    \label{dsd}
In any finite-dimensional representation of $N_+$, the matrix of $\ovl{x}$
satisfies:
\begin{equation}    \label{ds1}
\ovl{x}^{-1} (\pa_n + (\ovl{x} p_{-n} \ovl{x}^{-1})_-) \ovl{x} = \pa_n +
p_{-n} - \sum_{i \in I} c_i (p_{-n}^R \cdot \chi_i) p_i.
\end{equation}
\end{lem}

\begin{pf}
Using formula \eqref{actionr}, we obtain:
\begin{align*}
& \ovl{x}^{-1} (\pa_n + (\ovl{x} p_{-n} \ovl{x}^{-1})_-) \ovl{x} = \pa_n +
\ovl{x}^{-1} (p_{-n}^R \ovl{x}) + \ovl{x}^{-1} (\ovl{x} p_{-n}
\ovl{x}^{-1})_- \ovl{x} \\ & = \pa_n + \ovl{x}^{-1} (\ovl{x} p_{-n}
\ovl{x}^{-1})_+ \ovl{x} - \sum_{i \in I} c_i (p_{-n}^R \cdot \chi_i) p_i +
\ovl{x}^{-1} (\ovl{x} p_{-n} \ovl{x}^{-1})_- \ovl{x} \\ & = \pa_n + p_{-n}
- \sum_{i \in I} c_i (p_{-n}^R \cdot \chi_i) p_i,
\end{align*}
which coincides with \eqref{ds1}.
\end{pf}

\medskip
\noindent{\em Second proof of \thmref{ident}.} Let $K(z)$ be the map
${\Bbb A}^1 \arr N_+/A_+$ assigned to a smooth function ${\bold u}:
{\Bbb A}^1 \arr \h$ by \thmref{iso}. Let $\ovl{K}(z)$ be the element
of $N_+$ corresponding to $K(z)$ under the map $N_+/A_+ \arr N_+$
defined in \propref{localm}. By construction, $\ovl{K}(z)$ lies in the
$A_+$--coset of $K(z)$.

According to \lemref{p-1}, $(\ovl{K}(z) p_{-1} \ovl{K}(z)^{-1})_- =
p_{-1} + {\bold u}(z)$. Setting $n=1$ in formula \eqref{ds1} we
obtain:
$$\ovl{K}^{-1} (\pa_z + p_{-1} + {\bold u}(z)) \ovl{K} = \pa_z +
p_{-1} - \sum_{i \in I} c_i (p_{-1}^R \cdot \chi_i) p_i.$$ This shows
that $\ovl{K}(z)$ gives a solution to equation \eqref{ds}, and hence
lies in the $A_+$--coset of the Drinfeld-Sokolov dressing operator
$M(z)$. Therefore the cosets of $K(z)$ and $M(z)$ coincide.\qed
\medskip

It is possible to lift the map ${\bold u}(z) \arr N_+/A_+$ constructed
in \cite{DS} and \cite{FF} to a map ${\bold u}(z) \arr N_+$. We can
first attach to ${\bold u}(z)$ the coset $K(z)$ and then an element
$\ovl{K}(z)$ of $N_+$ defined as in the proof of \thmref{ident}. In
the next section we will show that $H_n = p_{-1}^R \cdot \chi_n \in
\C[u_i^{(n)}] \subset \C[N_+]$ (recall that $\chi_n
\not{\hspace*{-1mm}\in} \C[u_i^{(n)}]$). Since $p_{-1}^R \equiv
\pa_z$, we can view $\chi_n$ as $\int_{-\infty}^z H_n dz$. Hence we
can construct the image of ${\bold u}(z)$ in $N_+$ by the formula
$$\wt{K}(z) = \ovl{K}(z) \exp \left( \sum_{n\in I} c_n p_{n}
\int_{-\infty}^z H_n dz \right).$$ Comparing \eqref{ds1} and
\eqref{ds}, we can write an equivalent formula
$$M(z) \exp \left( - \sum_{n \in I} p_n \int_{-\infty}^z h_n(z) dz
\right),$$ where $h_n(z)$ are defined by formula \eqref{ds}. Note that the
last formula does not depend on the choice of $M(z)$.

We see that in contrast to the map to $N_+/A_+$, which is local,
i.e. depends only on the jet of ${\bold u}(z)$ at $z$, the map to
$N_+$ is non-local.

\begin{rem} In \cite{DS} it was proved that $h_n$ is the hamiltonian of the
$n$th equation of the mKdV hierarchy. In the next section we will prove in
a different way that $p_{-1}^R \cdot \chi_n$ is proportional to the
hamiltonian of the $n$th equation of the mKdV hierarchy.\qed
\end{rem}

\begin{rem} Now we can write an explicit formula for the Baker-Akhiezer
function associated to ${\bold u}$. Recall that this function is a
formal solution of equations \eqref{bakergen}. From formula
\eqref{ds1} we obtain the following solution:
\begin{align*}
\Psi({\bold t}) &= \ovl{K}({\bold t}) \exp \left( - \sum_{i \in I}
p_{-i}t_i - \sum_{i\in I} c_i p_i \int_{-\infty}^z H_i(z) dz \right)
\\ &= \wt{K}({\bold t}) \exp \left( - \sum_{i \in I} p_{-i}t_i \right),
\end{align*}
where ${\bold t} = \{ t_i \}_{i \in I}$ and $t_i$'s are the times of
the hierarchy (in particular, $t_1=z$). On the other hand, by
construction, the action of the vector field $\pa_n$ of the mKdV
hierarchy on $\wt{K} \in N_+$ corresponds to the right action of
$p_{-n} \in \ab_-$ on $N_+ \subset B_- \backslash G$. Hence if
$\wt{K}_0 \in N_+$ is the initial value of $\wt{K}$, when all $t_i=0$,
then $$\wt{K}({\bold t}) = \left( \wt{K}_0 \Gamma(\{ t_i \})
\right)_+,$$ where $$\Gamma({\bold t}) = \exp \left( \sum_{i \in I}
p_{-i}t_i \right)$$ and $g_+$ denotes the projection of $g \in B_-
\cdot N_+ \subset G$ on $N_+$ (it is well-defined for almost all
$t_i$'s). Finally, we obtain:
$$\Psi({\bold t}) = \left( \wt{K}_0 \Gamma({\bold t}) \right)_+
\Gamma({\bold t})^{-1}.$$

Note that this formula differs slightly from the one given in
\cite{Wil,HM} because in those papers another realization of the flag
manifold was chosen: $G/B_-$ instead of our $B_- \backslash G$.\qed
\end{rem}

\section{One-cocycles, hamiltonians and $\tau$-functions.}
In the previous section we established the equivalence between the
approaches of \cite{DS} and \cite{FF} to the mKdV hierarchies. In both
papers the mKdV equations were proved to be hamiltonian. In this
section we will discuss explicit formulas for the hamiltonians of the
mKdV equation and for some closely related cohomology classes of
$\n_+$. Note also that both in \cite{DS} and \cite{FF} it was shown
that the hamiltonians of the mKdV equations are integrals of motion of
the corresponding affine Toda equation (see also \cite{KW,W}).

\subsection{Connection between the hamiltonians and the
$\n_+$--cohomology.}  \label{conn}
In \cite{FF:laws,FF} the space spanned by the hamiltonians of the mKdV
equations was identified with the first cohomology of $\n_+$ with
coefficients in $\pi_0$, $H^1(\n_+,\pi_0)$. Let us briefly recall how to
assign an mKdV hamiltonian to a cohomology class.

The cohomology of $\n_+$ with coefficients in $\pi_0$ can be computed using
the Koszul complex $\pi_0 \otimes \bigwedge^*(\n_+^*)$. A cohomology class
from $H^1(\n_+,\pi_0)$ is represented in the Koszul complex by a functional
$f$ on $\n_+$ with coefficients in $\pi_0$, which satisfies the cocycle
condition $$f([a,b])-a\cdot f(b)+b\cdot f(a)=0.$$ This condition uniquely
determines $f$ by its values $f_i \in \pi_0$ on the generators $e_i,
i=0,\ldots,l$, of $\n_+$. Now set $g_i = \pa_z f_i - u_i f_i,
i=0,\ldots,l$. As shown in \cite{FF}, there exists $h \in \pi_0$, such that
$g_i = e_i^L \cdot h, i=0,\ldots,l$.

It was proved in \cite{FF:laws,FF} that $H^1(\n_+,\pi_0) \simeq
\ab_+^*$. Using the invariant inner product, we can identify $\ab_+^*$ with
$\ab_-$. Let $f_n$ be the cohomology class corresponding to $p_{-n} \in
\ab_-$. Then $h_n \in \pi_0$ constructed from $f_n$ is, by definition, the
density of the hamiltonian of the $n$th mKdV equation (i.e. the projection
of $h_n$ onto the space of local functionals $\pi_0/(\on{Im} \pa_z \oplus
\C)$ is an mKdV hamiltonian).

Below we give explicit formulas for $f_n(e_i)$ and $h_n$ as functions on
$N_+/A_+$. To simplify notation we will simply write $e_i$ for $e_i^L$ and
$p_{-n}$ for $p_{-n}^R$.

\subsection{Formulas for one-cocycles.}
Now recall that $\pi_0 \simeq \C[N_+/A_+]$. Hence the values of a
one-cocycle of $\n_+$ with coefficients in $\pi_0$ can be viewed as a
regular function on $N_+/A_+$.

\begin{prop}    \label{coc}
There exists a one-cocycle $\phi_n$ such that
\begin{equation}    \label{cocycle}
\left( \phi_n(e_i) \right)(K) = (e_i,Kp_{-n}K^{-1}), \quad \quad K
\in N_+/A_+.
\end{equation}
The cohomology classes corresponding to these cocycles span
$H^1(\n_+,\pi_0)$. In particular, if $n$ is a multiplicity free exponent,
then the cohomology classes defined by $\phi_n$ and $f_n$ coincide up to a
constant multiple.
\end{prop}

\begin{pf}
There exists a unique element $\De \in \HH \simeq \HH^*$,
such that $(\al_i,\De)=1,
\forall i=0,\ldots,l$, and $(d,\De)=0$. But $\HH^*$ is isomorphic to $\HH$
via the non-degenerate inner product $(\cdot,\cdot)$. Let us use the same
notation for the image of $\De$ in $\HH$ under this isomorphism. Then $\De$
satisfies: $[\De,e_i]=e_i, [\De,h_i]=0, [\De,f_i] = -f_i,
i=0,\ldots,l$. Thus, the adjoint action of $\De$ on $\G$ coincides with the
action of the principal gradation.

Any function $F \in \C[N_+]$ can be viewed as an element of the zeroth
group of the Koszul complex of the cohomology of $\n_+$ with coefficients
in $\C[N_+]$. The coboundary of this element is a (trivial) one-cocycle,
whose value on $e_i$ is $e_i \cdot F \in \C[N_+], i=0,\ldots,l$.

Consider a function $F_n=\De_n$ on $N_+$ defined by the formula
\begin{equation}    \label{Fn}
F_n(\xx )=(\De,\xx  p_{-n}\xx^{-1})
\end{equation}
Note that here $p_{-n}$ is considered as an element of $\G$ and we consider
the adjoint action of $N_+$ on $\G$. The value of the corresponding
one-cocycle on $e_i$ is equal to $e_i \cdot F_n$. We have: $$(e_i \cdot
F_n)(\xx ) = (\De,[e_i,\xx p_{-n} \xx^{-1}]) = ([\De,e_i],\xx p_{-n}
\xx^{-1}) = (e_i,\xx p_{-n} \xx^{-1}).$$ Thus, there exists a one-cocycle
$f$ of $\n_+$ with coefficients in $\C[N_+]$, such that
\begin{equation}    \label{esch}
f(e_i) = (e_i,\xx p_{-n} \xx^{-1}), i=0,\ldots,l.
\end{equation}
Moreover, $f(e_i)$ is $A_+$--invariant for all $i=0,\ldots,l$. Indeed,
$$(p_m \cdot f(e_i))(x) = (e_i,\xx [p_m,p_{-n}] \xx^{-1}) = n(p_n,p_{-n})
(e_i,\xx C \xx^{-1}) \delta_{n,-m}$$ $$= n(p_n,p_{-n}) (e_i,C)
\delta_{n,-m} =0.$$ Therefore formula \eqref{esch} defines a one-cocycle of
$\n_+$ with coefficients in $\C[N_+/A_+] \simeq \pi_0$.  This is the
cocycle $\phi_n$.

By construction, $\phi_n$ is a trivial one-cocycle of $\n_+$ with
coefficients in $\C[N_+]$. But it is non-trivial as a one-cocycle of
$\n_+$ with coefficients in $\C[N_+/A_+]$. Indeed, if it were a
coboundary, there would exist an $A_+$--invariant function $\wt{F}_n$
on $N_+$, such that $\phi_n(e_i)=e_i \cdot \wt{F}_n$. But then $e_i
\cdot (\wt{F}_n-F_n)=0$ for all $i$, and $\wt{F}_n-F_n$ is
$N_+$--invariant, and hence constant. However, by \eqref{cox}, $p_n
\cdot F_n = nh (p_n,p_{-n}) \neq 0$, where $h$ is the Coxeter number
of $\G$. Hence the function $F_n$ is not $A_+$--invariant.

Thus, $\wt{F}_n-F_n$ can not be a constant function. Therefore $\phi_n$
defines a non-zero cohomology class. Let us compute its degree with respect
to the principal gradation. We have:
$$\left( \De \cdot (\phi_n(e_i)) \right)(\xx ) = (e_i,[(\xx \De
\xx^{-1})_+,\xx p_{-n}\xx^{-1}])$$ $$= (e_i,[\xx \De \xx^{-1},\xx
p_{-n}\xx^{-1}]) - (e_i,[(\xx \De \xx^{-1})_-,\xx p_{-n}\xx^{-1}])$$
$$=(e_i,\xx [\De,p_{-n}]\xx^{-1}) - ([e_i,\De],\xx p_{-n}\xx^{-1}]) =
(-n+1) \phi_n(e_i).$$ Hence the degree of $\phi_n$ equals $-n$.

For multiplicity free exponent $n$ this implies that the cohomology class
of $\phi_n$ is proportional to that of $f_n$. For the multiple exponents
$i$ which occur in the case of $D_{2n}^{(1)}$ (see \remref{multi}), we need
to show that the cocycles $\phi_i^1$ and $\phi_i^2$, corresponding to two
linearly independent elements $p_{-i}^1$ and $p_{-i}^2$ of $\ab_-$ of
degree $-i$, are linearly independent. But a linear combination $\alpha
\phi_i^1+\beta \phi_i^2$ of these cocycles is just the cocycle
corresponding to $\alpha p_{-i}^1+\beta p_{-i}^2 \in \ab_-$. The argument
that we used above can be applied to the cocycle $\alpha \phi_i^1+\beta
\phi_i^2$ to show that it is non-trivial unless both $\alpha$ and $\beta$
equal $0$.
\end{pf}

\begin{rem} Homogeneous functions on $N_+/A_+$ are necessarily
algebraic. Thus, $\phi_n(e_i)$ $\in$ $\C[N_+/A_+]$.\qed
\end{rem}

\begin{rem} One can show in the same way as above that for any $\chi \in
\HH$, such that $(\chi,C) \neq 0$, there exists a one-cocycle
$\wt{\chi}_n$ of $\n_+$ with coefficients in $\C[N_+/A_+]$, which
satisfies the following property: considered as an $A_+$--invariant
function on $N_+$, $\wt{\chi}_n(e_i)$ equals $e_i \cdot \chi_n$, where
$\chi_n$ is the function on $N_+$ defined in Sect.~4.1. The
one-cocycle $\wt{\chi}_n$ is homologous to $F_n$, suitably
normalized.\qed
\end{rem}

\subsection{Formulas for hamiltonians.}
Now we can find a formula for the density of the $n$th mKdV
hamiltonian using \propref{coc} and the procedure of \secref{conn}.

\begin{prop} The function $H_n$ on $N_+/A_+$, such that
$$H_n(K) = (p_{-1},Kp_{-n}K^{-1}), \quad \quad K \in N_+/A_+$$ is a
density of the $n$th hamiltonian of the mKdV hierarchy.
\end{prop}

\begin{pf} We have to show that
\begin{equation}    \label{prel}
e_i \cdot H_n = p_{-1} \phi_n(e_i) - u_i \phi_n(e_i), \quad \quad
i=0,\ldots,l.
\end{equation}

Let us consider functions on $N_+/A_+$ as $A_+$--invariant functions on
$N_+$. Recall from Sect.~2 that there is a unique up to a constant
isomorphism $\epsilon_\la$ between $\C[N_+]$ and the contragradient Verma
module $M_{\la}^*$, which commutes with the left action of $\n_+$. For $a
\in \g$ the operator $\ep_\la a \ep_\la^{-1}$ on $\C[N_+]$ is the first
order differential operator $a^R+f_\la(a)$. Here $f_\la(a) = \ep_\la^{-1}
(a \cdot v_\la)$ We know that $u_i = f_{\al_i}(p_{-1})$, see \cite{FF} and
Sect.~2. Hence $\ep_{-\al_i}^{-1} p_{-1} \ep_{-\al_i} = p_{-1}^R-u_i$, and
hence formula \eqref{prel} can be rewritten as
\begin{equation}    \label{cob}
\ep_{-\al_i} (e_i \cdot H_n) = p_{-1} \cdot \ep_{-\al_i}(\phi_n(e_i)) \quad
\quad i=0,\ldots,l.
\end{equation}

Let us show that $H_n = p_{-1} \cdot F_n$, where the function $F_n \in
\C[N_+]$ is defined by formula \eqref{Fn}. Indeed, $$(p_{-1} \cdot F_n)(x)
= (\De,[(xp_{-1}x^{-1})_+,xp_{-n}x^{-1}]) = -
(\De,[(xp_{-1}x^{-1})_-,xp_{-n}x^{-1}])$$ $$= -
([\De,(xp_{-1}x^{-1})_-],xp_{-n}x^{-1}) = ([\De,p_{-1}],xp_{-n}x^{-1}) =
(p_{-1},xp_{-n}x^{-1}).$$

The fact that $H_n$ is $A_+$--invariant can be proved in the same way as
for $\phi_n(e_i)$.

Now recall that the map $\ep_{-\al_i} e_i: \C[N_+] \arr M_{-\al_i}^*$
commutes with the action of $\g$, where $\g$ acts on $\C[N_+]$ from
the right by vector fields, see \cite{FF}, Sect.~4.  Therefore we
obtain
$$\ep_{-\al_i} e_i (p_{-1} \cdot F_n) = p_{-1} \cdot \ep_{-\al_i} e_i
(F_n).$$ This implies formula \eqref{cob} if we take into account that $H_n
= p_{-1} \cdot F_n$ and $\phi_n(e_i) = e_i \cdot F_n$.
\end{pf}

\begin{rem} Our formula for the hamiltonians is
equivalent to the formula given by Wilson \cite{W}, (4.10).\qed
\end{rem}

\begin{rem} One can also construct the density of the $n$th hamiltonian
as $p_{-1} \cdot \chi_n$ where $\chi_n$ was defined in Sect.~4.1. For
different $\chi$, these densities, suitably normalized, differ by
total derivatives, and hence define the same hamiltonian.\qed
\end{rem}

\subsection{Involutivity of the hamiltonians.}
Now we want to prove that the Poisson bracket between two mKdV hamiltonians
vanishes. This is equivalent to showing that $p_{-n} \cdot H_m = p_{-1}
H_{n,m}$ for some $H_{n,m} \in \C[N_+/A_+]$ (see \cite{FF}).

\begin{prop} Define $H_{n,m} \in \C[N_+/A_+]$ by formula $$H_{n,m}(K) =
-([\De,(Kp_{-n}K^{-1})_-],Kp_{-m}K^{-1}).$$ Then $$p_{-n} \cdot H_m = p_{-m}
\cdot H_n = p_{-1} \cdot H_{n,m}.$$
\end{prop}

\begin{pf} We have: $$(p_{-n} \cdot H_m)(K) =
(p_{-1},[(Kp_{-n}K^{-1})_+,Kp_{-m}K^{-1}]) =$$ $$
(p_{-1},[(Kp_{-n}K^{-1})_+,(Kp_{-m}K^{-1})_-]),$$ because
$(p_{-1},[y_1,y_2])=0$ if $y_1,y_2 \in \n_+$. On the other hand, $$(p_{-m}
\cdot H_n)(K) = (p_{-1},[(Kp_{-m}K^{-1})_+,Kp_{-n}K^{-1}]) =$$ $$-
(p_{-1},[(Kp_{-m}K^{-1})_-,Kp_{-n}K^{-1}]) = -
(p_{-1},[(Kp_{-m}K^{-1})_-,(Kp_{-n}K^{-1})_+]),$$ because $(p_{-1},y)=0$ if
$y \in \bb_-$. Therefore $p_{-n} \cdot H_m = p_{-m} \cdot H_n$.

Consider now $H_m$ as an $A_+$--invariant function on $N_+$. Then we have:
$H_m = p_{-1} \cdot F_m$. Hence $p_{-n} \cdot H_m = p_{-1} \cdot (p_{-n}
\cdot F_m)$. Let $H_{n,m} = p_{-n} \cdot F_m$. We obtain: $$(p_{-n} \cdot
F_m)(x) = (\De,[(xp_{-n}x^{-1})_+,xp_{-m}x^{-1}]) = -
(\De,[(xp_{-n}x^{-1})_-,xp_{-m}x^{-1}])$$ $$= -
([\De,(xp_{-n}x^{-1})_+],xp_{-m}x^{-1}).$$ The latter expression is
$A_+$--invariant, which can be shown in the same way as in the proof of
\propref{coc}. Hence $H_{n,m} \in \C[N_+/A_+]$ and $p_{-n} \cdot H_m =
p_{-m} \cdot H_n = p_{-1} H_{n,m}$.
\end{pf}

\subsection{Connection with $\tau$--functions.}
The $\tau$--functions have the following meaning from our point of
view. For $\la \in \wt{\h}^*$, consider the contragradient Verma
module $M_\la^*$ over $\g$. This module can be realized in the space
of sections of a line bundle $\xi_\la$ over $N_+$, considered as a
big cell of the flag manifold $\fl$. By definition, the
$\tau$--function $\tau_\la$ corresponding to $\la$ is the unique up
to a constant $N_+$--invariant section of $\xi_\la$ over $N_+$.

\begin{rem} This should be compared with the definition of the
$\tau$--functions in the framework of the Grassmannian approach
\cite{DJKM,SW,Wil}.\qed
\end{rem}

Note that $\xi_\la$ can be trivialized over $N_+$, and so there exists a
unique up to a non-zero constant isomorphism between the space of sections
of $\xi_\la$ and $\C[N_+]$. Under this isomorphism, $\tau_\la$
corresponds to a constant function on $N_+$.

Let $\La_i, i=0,\ldots,l$, be the fundamental weights of the affine algebra
$\g$. We call $\tau_{\La_i}$ the $i$th $\tau$--function of $\g$ and denote
it by $\tau_i$. Let us also set $\tau = \tau_{\De}$.

According to \propref{kost}, for any $a \in \g$, $a \cdot
\tau_\la = f_\la(a) \tau_\la$, where $f_\la(a)(\xx ) = \langle
\la,\xx a\xx^{-1} \rangle$.

In particular, we see that $e^{\varphi_i}=\tau_{\al_i}$, and $p_{-1}
e^{\varphi_i} = \pa_z e^{\varphi_i} = u_i e^{\varphi_i}$. Note that
$e^{\varphi_i}$ can be expressed in terms of $\tau_j$'s. For example, for
$\g=\widehat{\sw}_N$ we have: $e^{\varphi_i} = \tau_{i-1}^{-1} \tau_i^2
\tau_{i+1}^{-1}$, which is well-known.

Now we can interpret the functions $F_n$ as logarithmic derivatives of
$\tau$. Indeed, we obtain $\pa_n \tau = F_n \tau$, so that we can formally
write: $F_n = \pa_n \log \tau$. Further, $H_n = \pa_n \pa_z \log
\tau$, and, more generally, $H_{n,m} = \pa_n \pa_m \log \tau$, which
coincides with known results. Similarly, we can write: $u_i^{(n)} =
\pa_z^{n+1} \log \tau_{\al_i}$.

\begin{rem} More generally, we have the following formula for the function
$\chi_n$ defined in Sect.~4.1: $\chi_n = \pa_n
\tau_\chi/\tau_\chi$.\qed
\end{rem}

To summarize, the group $N_+$ has natural coordinates $u_i^{(n)}$ and
$F_n$, which can be obtained as logarithmic derivatives of
$\tau$--functions. The vector fields $p^R_{-n}$ written in terms of
these coordinates provide the flows of the mKdV hierarchy, and the
vector field $\sum_{i=0}^l e^L_i$ written in terms of these
coordinates gives the affine Toda equation \cite{FF}.

\subsection{Example of $\su$.} Here we will write explicit formulas for
the action of the generators of the nilpotent subalgebra of $\su$ and mKdV
hamiltonians on the corresponding unipotent subgroup.

According to the results of this section, we have an isomorphism $$\C[N_+]
\simeq \C[u^{(n)},F_m]_{n\geq 0,m \on{odd}}.$$ The left action
of the generators $e_0$ and $e_1$ of $\n_+$ on $\C[N_+]$ is given by

\begin{align*}
e_0 &= - \sum_{n\geq 0} P^+_n \frac{\pa}{\pa u^{(n)}} + \sum_{m \on{odd}}
\phi_m(e_0) \frac{\pa}{\pa F_m},\\
e_1 &= - \sum_{n\geq 0} P^-_n \frac{\pa}{\pa u^{(n)}} + \sum_{m \on{odd}}
\phi_m(e_1) \frac{\pa}{\pa F_m},
\end{align*}
where $P^\pm_n$ are elements of $\C[u^{(n)}]$, defined recursively as
follows: $P^\pm_0 = 1$, $P^\pm_{n+1} = \pa P^\pm_n \pm u P^\pm_n$, and
$\phi_m(e_i)$ are the values of a one-cocycle $\phi_m$ of $\n_+$ with
coefficients in $\C[u_i^{(n)}]$ of degree $m$.

The right action of $p_k, k$ positive odd, is given by $4k \pa/\pa F_k$,
and the action of $p_{-k}, k$ positive odd, is given by
$$p_{-k} = \sum_{n\geq 0} (\pa^{n+1} q_k) \frac{\pa}{\pa u^{(n)}} +
\sum_{m \on{odd}} H_{k,m} \frac{\pa}{\pa F_m}.$$

The $m$th mKdV equation now reads: $$\pa_m u = q_m.$$ This equation is
hamiltonian with the hamiltonian $(1/m) H_{m,1}$, and hence
$$q_m = \frac{1}{m} \frac{\delta H_{m,1}}{\delta u}.$$

The involutivity of the hamiltonians means that $$\sum_{n\geq 0} (\pa^{n+1}
q_k) \frac{\pa H_{m,1}}{\pa u^{(n)}} = \pa H_{k,m}$$ (note that $H_{k,m} =
H_{m,k}$ and $H_{m,1} = H_{1,m} = H_m$).

The KdV variable is $v=\frac{1}{2}u^2+u'$, and $\C[v^{(n)}]_{n\geq 0}
\subset \C[u^{(n)}]_{n\geq 0}$ coincides with the $e_1$--invariant
subspace of $\C[u^{(n)}]_{n\geq 0}$.

\noindent{\bf Acknowledgements.} We would like to thank B.~Feigin for
his collaboration in \cite{EF,FF:laws,FF} and useful discussions. The
second author also thanks D.~Ben-Zvi for interesting discussions.

The research of the second author was supported by grants from the
Packard Foundation, NSF and the Sloan Foundation.

\end{document}